\documentstyle[11pt, a4wide, leqno]{article}

\def\bea{\begin{eqnarray}}
\def\eea{\end{eqnarray}}
\def\beann{\begin{eqnarray*}}
\def\eeann{\end{eqnarray*}}
\def\beq{\begin{equation}}
\def\eeq{\end{equation}}
\def\ba{\begin{array}}
\def\ea{\end{array}}
\def\ben{\begin{enumerate}}
\def\een{\end{enumerate}}

\newtheorem{th}{Theorem}[section]

\newtheorem{pro}[th]{Proposition}

\newtheorem{co}[th]{Corollary}

\font\mybb=msbm10 at 11pt

\def\bb#1{\hbox{\mybb#1}}

\def\bZ {\bb{Z}}
\def\bR {\bb{R}}

\def\bC {\bb{C}}

\def\dd{d^\dagger}

\date{}

\begin{document}
\rightline{CERN-TH/2000-282}

\thispagestyle{empty}
\begin{center}
\vspace {.7cm} {\large {\bf{ Vanishing Theorems and String
Backgrounds}} } \vskip  1.0truecm {\large{\bf{ S. Ivanov${}^1$ and
G. Papadopoulos${}^2$}}} \vskip 2.0truecm
 {\normalsize{\sl 1. Department of Mathematics}}
\\ {\normalsize{\sl University of Sofia}}
\\ {\normalsize{\sl \lq\lq St. Kl. Ohridski''}}
\vskip 1.0truecm
{\normalsize{\sl 2. CERN}}
\\ {\normalsize{\sl Theory Division}}
\\ {\normalsize{\sl 1211 Geneva, 23}}
\\ {\normalsize{\sl Switzerland}}

\vskip 2.0truecm
{\bf {Abstract}}
\end{center}
We show various vanishing theorems for the cohomology groups  of
compact hermitian manifolds for which the Bismut connection has
(restricted) holonomy contained in $SU(n)$ and   classify  all
such manifolds of dimension four. In this way we provide necessary
conditions for the existence of such structures on hermitian
manifolds.  Then we apply our results to solutions of the string
equations and show that such solutions admit various cohomological
restrictions like for example that under certain natural
assumptions  the plurigenera vanish. We also find that under some
assumptions the string equations are equivalent to the condition
that a certain vector is parallel with respect to the Bismut
connection.

\newpage

\section{Introduction}

Riemannian manifolds equipped with a closed form have found many
applications in various branches of mathematics and physics. In
physics, the classical example is that of manifolds equipped with
a closed two-form which describe gravity in the presence of a Maxwell
field. More recently, Riemannian or pseudo-Riemannian manifolds
$M$ equipped  with (closed) forms of any degree arise in the
context of string theory. In particular in string
compactifications, one investigates manifolds $\bR^k\times M$,
where $M$ is a compact manifold of dimension $10-k$. For most
applications, the closed forms   vanish and $M$ is required to
have special holonomy \cite{candelas, gpt}. Then  much  of the
physical data, like the  particle spectrum on $\bR^k$, which are
extracted from these compactifications, are determined in terms of
the cohomological information on $M$, like for example Euler and
Hodge numbers.  Recently  however, compactifications have also
been considered for which the closed forms do not vanish
\cite{strom, becker, vafa}. Amongst the various closed forms that
appear in string theory, there is a
 closed three-form $H$ which is associated
to the fundamental string. Such compactifications
have been investigated in \cite{strom}. In
this case,  this three-form can be
interpreted as torsion of a
metric connection, $\nabla$, on $M$.
For compactifications of string theory with non-vanishing
$H$, it is required that the  connection with torsion $\nabla$ has
  holonomy which is contained
in $SU(n)$ ($n=2,3,4$), $G_2$ or  ${\rm Spin}(7)$.
The case of  holonomy $SU(n)$ is of particular interest
because the underlying manifolds $M$  are  hermitian.
The connection with torsion can then be identified
with the Bismut connection \cite{Bi}; for the various
 connections  on hermitian
manifolds see \cite{yano,yanko}. Other applications include the
geometry on the target spaces of supersymmetric one- \cite{coles}
and two-dimensional sigma models \cite{rocek, hp1}. More recently
it has been shown that the geometry on the moduli space of certain
five-dimensional black holes  is hyper-K\"ahler with torsion
(HKT)\cite{HP2}, ie the  holonomy of the Bismut
connection is in $Sp(k)$ \cite{GPS, strom2, jggp}, .

In mathematics, a natural arena for the
investigation of metric connections
 with torsion a three-form is the theory of Hermitian manifolds.
Hermitian manifolds apart from the Bismut connection that has been
mentioned above also admit the Chern connection (see for example
\cite{yano,yanko}). Both connection have  holonomy
which is contained in $U(n)$. The Chern connection, in addition,
has curvature which is a (1,1)-form with respect to the natural
complex structure and therefore it is compatible with the
holomorphic structure of the tangent bundle. Consequently, many
theorems regarding the non-existence of holomorphic sections have
been expressed in terms of conditions on the curvature of the
Chern connection. Recently, however, vanishing theorems for
Dolbeault cohomology groups have been expressed in terms of
conditions  on the curvature of the Bismut connection \cite{AI1}.

In string theory, the metric $g$
and the three-form $H$ of a manifold $M$ are required
to satisfy the \lq string equations'
which are a generalization
of the Einstein equations of general relativity. For this, in
addition,  a function
$\phi$ is introduced on $M$, called dilaton, and the (type II)
 string equations
can be written as follows:
\begin{eqnarray}\label{st}
& &Ric^g_{ij} -{1\over4} H_{imn}H_j^{mn} + 2\nabla^g_i\partial_j\phi =0,\\
& &\nabla^g_i\left(e^{-2\phi}H^{imn}\right)=0\ ,\nonumber
\end{eqnarray}
where $Ric^g$ and $\nabla^g$ are the Ricci
tensor and the Levi-Civita connection
of the metric $g$. The three-form $H$ is  closed
in the context of (type II) string theory but we can take it to be any
three-form. In what follows, if $H$ is required to be closed, it will
be mentioned explicitly. In fact there is a third field equation
that associated with the dilaton $\phi$. However the dilaton
equation is implied from those of (\ref{st}) up to a
constant and so we shall not further investigate it.
The above equations are also supplemented with
two additional equations, called Killing spinor
equations,  the following:
\begin{eqnarray}\label{kse}
\nabla\eta&=&0
\nonumber\\
(d\phi-{1\over6}H)\eta&=&0\ ,
\end{eqnarray}
where $\eta$ is a spinor and we have denoted with the same symbol
the forms and their associated Clifford algebra element. In fact
in type II string theory there is an additional set of Killing
spinor equations for another spinor $\epsilon$ which are related
to those in (\ref{kse}) by taking $H$ to $-H$. The solutions of
the string equations that are of interest are those for which
there is a non-vanishing spinor $\eta$ which satisfies the Killing
spinor equations. In the heterotic string theory, the situation is
somewhat different. The form $H$ is not necessarily closed.
In fact the \lq sigma model anomaly cancellation at one loop' requires
that $dH$  is proportional to the difference
of the first Pontriangin classes of $M$ and a vector bundle
over $M$ with coefficient which depends on an expansion parameter
called  string tension. Apart from
the Killing spinor equations of the heterotic string   given in
(\ref{kse}), there is also another one associated with the gauge
sector of the theory which however does not affect most of our
considerations. It
will be again mentioned  though in 
sections four and  five. 
In both type II
and heterotic strings, the string and killing spinor
equations above receive higher order corrections which typically 
involve powers of the curvature tensor with expansion parameter
the string tension. In fact the string equations for the
heterotic string involve additional terms arising from
the gauge sector even in the same order as that of (\ref{st}).

Suppose that $M$ is a KT manifold. The first Killing spinor
equation in (\ref{kse}) requires that the spinor is parallel with
respect to the Bismut connection. Manifolds that admit such
parallel spinors are those  for which the  holonomy of
the Bismut connection is a subgroup of $SU(n)$; for applications
in  string theory $n\leq 4$; for n=2 the existence of parallel
spinors with respect to the Bismut connection leads to some
constraints on the KT geometry depending also on the type of the
parallel spinor \cite{daliv}. The second Killing spinor equation
in (\ref{kse}) imposes additional conditions which have been
investigated in \cite{strom}. Here we shall not be concerned with
solutions of this equation; we have already given some
applications of this in \cite{SIGP}. In any case the second
Killing spinor equation does not directly arise in the
world-volume (conformal field theory) approach to strings. This is
unlike the first one which is a sufficient condition for the
existence of a complex structure which characterizes the relevant
conformal field theories.

In the first part of the paper, we shall first establish various
relations between the curvature tensors of the Chern and Bismut
connections. Then  we shall state various vanishing theorems, like
that of Kodaira,  for hermitian manifolds in terms of conditions
on the curvature of the Bismut connection. We shall find that if
the (restricted) holonomy of the Bismut connection is contained in
$SU(n)$ and a certain condition on the torsion is satisfied, then
the Kodaira-type vanishing theorem for holomorphic (p,0)-forms
holds. As an application of our results we show that a certain
class of balanced hermitian manifolds does not admit such KT
structures. Then we shall classify all four-dimensional compact
hermitian manifolds for which the Bismut connection has
(restricted) holonomy contained in $SU(2)$. We shall show that
such spaces are either conformal to a Calabi-Yau space or to  a
Hopf surface.

In the second part of the paper, we shall apply the above
results to investigate the cohomological properties of the
solutions of the string equations. We shall first find that
under the assumption of $SU(n)$ holonomy for the Bismut
connection, the string equations can be expressed
in a simple form. Then  we shall show that some
of the cohomology groups of such spaces vanish.

This paper is organized as follows: In section two, we give
various definitions of the manifolds that we shall investigate and
establish our notation. In section three, we show various
identities that relate the curvature of the Chern and Bismut
connections. In section four, we present our main theorems and
apply our results to (i) balanced hermitian manifolds and (ii)
four-dimensional hermitian manifolds which admit a Bismut
connection with (restricted) holonomy contained in $SU(2)$. In
section five, we investigate various compact solutions to the
string equations and describe certain properties of their
cohomology groups. In section six, we give our conclusions.

\section{Hermitian and KT Manifolds}

Let $(M,g,J)$ be a $2n$-dimensional ($n>1$) Hermitian manifold with complex
structure $J$ and compatible Riemannian metric $g$. The
K\"ahler form $\Omega$ of $(M,g,J)$ is defined by
\beq
\Omega (X,Y)=g(X,JY)\ .
\eeq
 Denote by $\theta$ the Lee form of $(M,g,J)$,
\beq \theta = \dd \Omega \circ J\ , \eeq where $\dd$ is the
adjoint of $d$. For a one-form $\alpha$, we shall denote by
$J\alpha $ the form dual to $J\alpha^{\#}$, where $\alpha ^{\#}$
is the vector dual to $\alpha $. Equivalently, $J\alpha = -\alpha
\circ J$. Hence, $\dd \Omega = J\theta $. It has been shown by
Gauduchon \cite{G3} that  any conformal class of Hermitian metrics
on a compact manifold contains a unique (up to homothety) metric
satisfying $\dd \theta =0$. This metric is called the {\it
Gauduchon metric}.

The Bismut connection $\nabla$ and the Chern connection $D$
are given by
\begin{equation}\label{1}
g(\nabla_X Y,Z) = g(\nabla^g_X Y,Z) + \frac{1}{2} d^c \Omega (X,Y,Z),
\end{equation}
\begin{equation}\label{1dop}
g(D_X Y,Z) = g(\nabla^g_X Y,Z) + \frac{1}{2} d \Omega (JX,Y,Z),
\end{equation}
respectively, where $\nabla^g$ is the Levi-Civita connection of $g$.
Recall that $d^c = i(\overline {\partial} - \partial )$.
In particular,
$d^c \Omega (X,Y,Z) = - d \Omega (JX,JY,JZ)$. Both these
connections on Hermitian manifolds have been known for
sometime, see for example \cite{yano, yanko, Bi, Ga1}.

Let $T$ be the torsion of $\nabla$,
$T(X,Y)=\nabla_XY-\nabla_YX-[X,Y]$, and $C$ be the torsion of $D$,
respectively. It follows from the definition of the Bismut
connection (\ref{1})  that \bea \nabla g& = &0, \quad \nabla J=0,
\nonumber\\ T(X,Y,Z):& = &g(T(X,Y),Z)=d^c\Omega (X,Y,Z);
\label{1a} \eea The equation (\ref{1a}) shows that $T=d^c\Omega$,
$dT=dd^c\Omega=2i
\partial\bar{\partial}\Omega$. Hermitian manifolds equipped
with the Bismut connection are called {\it K\"ahler with torsion} (KT).
A hermitian manifold admits a {\it strong} KT structure if the
torsion of the Bismut connection is {\it closed}. Hence, a
 KT manifold is strong, iff its K\"ahler
form is $\partial\bar{\partial}$-closed.

{}From the definition  (\ref{1dop}) of the Chern connection, we have
\bea\label{2dop}
Dg& = &0, \quad DJ=0,
\nonumber\\
 2C(X,Y,Z)
:& = & 2g(C(X,Y),Z)=d\Omega(JX,Y,Z)+d\Omega(X,JY,Z).
\eea
The equality (\ref{2dop}) yields $C(JX,Y)=C(X,JY)$ which implies \cite{Bal}
 $C(JX,Y)=JC(X,Y)$.

Using the formula (see e.g. \cite{KN}) $$ (\nabla^g_X\Omega)(Y,Z)
= -g\left((\nabla^g_XJ)Y,Z\right) = - \frac{1}{2}
\left(d\Omega(X,JY,JZ)-d\Omega(X,Y,Z)\right)\ , $$
 we get the expressions
\begin{equation}\label{3}
\theta (X)=\dd\Omega(JX) = -\frac{1}{2}\sum _{i=1}^{2n}T(JX,e_i,Je_i) =
\frac{1}{2}\sum _{i=1}^{2n}C(JX,e_i,Je_i).
\end{equation}
Here and henceforth $e_1,e_2,...,e_{2n}$
is an orthonormal basis of the
tangential space.

Let $(M,g,(J_a), a=1,2,3)$ be a 4n-dimensional hyper-K\"ahler manifold
with torsion  (HKT) \cite{HP2}.
We denote with $\Omega_a$ and $\theta_a$ the  K\"ahler form and the Lee form
of the complex structure $J_a$, respectively.
For HKT manifolds,   the
Bismut connections of the three
hermitian structures $(g,J_a), a=1,2,3$ associated
with the hypercomplex structure $J_a, a=1,2,3$ coincide and this
condition
is equivalent to \cite{HP2,GPS,Ga1,GrP}
$$
d_1\Omega_1=d_2\Omega_2=d_3\Omega_3\ ,
$$
where $d_a\Omega_a(X,Y,Z):=-
d\Omega_a(J_aX,J_aY,J_aZ), a=1,2,3$.
In particular, if one of the hermitian
structures  $(g,J_a), a=1,2,3$ of $M$ is K\"ahler,
then the other two are also K\"ahler
and so $M$ is a hyper-K\"ahler manifold.
The torsion $T$ is $(2,1)+(1,2)$-form
with respect to each complex
structure $J_a, a=1,2,3$. This leads to the
equality  of the three  Lee forms
\cite{GT,I1}, i.e. $\theta_1=\theta_2=
\theta_3:=\theta$.

Let $R=[\nabla,\nabla]-\nabla_{[,]}$ be
the curvature tensor of type (1,3) of
the Bismut connection $\nabla$. We
denote the curvature tensor of type (0,4)
 $R(X,Y,Z,V)=g(R(X,Y)Z,V)$ by
the same letter. We shall denote the curvature of the Chern
and the Levi-Civita
connections with $K$ and $R^g$, respectively.

The Ricci tensor $Ric$ and the  Ricci form $\rho$ of the Bismut
connection $\nabla$ are defined  by
$$
Ric(X,Y) = \sum_{i=1}^{2n}R(e_i,X,Y,e_i), \quad
\rho(X,Y)=\frac{1}{2} \sum_{i=1}^{2n}R(X,Y,e_i,Je_i)\ .
$$

The two Ricci forms $\rho^D$ and $\kappa$ associated
with the Chern connection are defined by
$$
\rho^D(X,Y)=\frac{1}{2} \sum_{i=1}^{2n}K(X,Y,e_i,Je_i), \quad
\kappa(X,Y)=\frac{1}{2} \sum_{i=1}^{2n}K(e_i,Je_i,X,Y).
$$

The (1,1)-form $\rho^D$ represents the first Chern class of the
manifold $M$.
The (1,1)-form $\kappa$ is called sometimes \lq the mean curvature' of
the holomorphic tangent bundle $T^{1,0}M$ with the  hermitian
metric induced by $g$ \cite{Ko1}. Vanishing theorems for holomorphic
(p,0)-forms $p=1,2,\ldots,n$, which we shall use later, are expressed in
terms of the non-negativity of $\kappa$ \cite{KW,Ga}.

The two Ricci forms $\rho$ and $\rho^D$ are related by \cite{AI1}
\begin{equation}\label{ricc1}
\rho^D =\rho + d(J\theta).
\end{equation}
This can be most easily proved by computing first the trace of the
Chern and Bismut connections with $J$ and then comparing the
two expressions using the definition of $\theta$.

We shall adopt the conventions of \cite{AI1} to denote the trace
of $\rho$ with $J$ by $b$ and the conventions of \cite{G4} to
denote   the trace of $\kappa$ with $J$, which is  equal to the
trace of $\rho ^D$, by $2u$, ie \bea\label{sca}
b&=&\sum_{j=1}^{2n} \rho (Je_j,e_j), \nonumber\\
2u&=&\sum_{j=1}^{2n} \rho ^D (Je_j,e_j)=\sum_{j=1}^{2n}
\kappa(Je_j,e_j)\ . \eea The standard scalar curvature
$Scal^{\nabla}$ of the Bismut connection $\nabla$, is defined by
$$ Scal^{\nabla}=\sum_{j=1}^{2n} Ric (e_j,e_j)\ . $$

The trace of the exterior derivative
$dT$ of the torsion 3-form with $J$,
 we denote
by $\lambda^\Omega$, i.e.
\begin{equation}\label{tr}
\lambda^{\Omega}(X,Y) = \sum_{i=1}^{2n} dT(X,Y,e_i,Je_i).
\end{equation}
Note that $\lambda^{\Omega}$ is an
(1,1)-form with respect to $J$ since $dT$ is a
(2,2)-form.

In what follows the following definition seem to be useful\\[1mm]

{\bf Definition} A KT (resp. HKT) manifold is
said to be {\it almost strong KT}
(resp. {\it almost strong HKT}) manifold if $\lambda^\Omega=0$
(resp.  $\lambda^{\Omega_1}=0$).
\\[1mm]

Note that for  HKT manifolds
$\lambda^{\Omega_1}(J_1.,.)=\lambda^{\Omega_2}(J_2.,.)=
\lambda^{\Omega_3}(J_3.,.)$
since the four-form $dT$ is of type (2,2) with respect
 to each of the three complex structures (see e.g. \cite{I1}).
Clearly, every strong KT manifold is almost strong KT.
If a 2n-dimensional KT manifold is locally conformally K\"ahler, then
$T=\frac{1}{n-1}J\theta \wedge \Omega$. In such case, a straight
forward computation
reveals that
\begin{equation}\label{four1}
(n-1)\lambda^\Omega = (4-2n)(dJ\theta+\theta \wedge J\theta +
|\theta|^2\Omega)-2\dd\theta \Omega,
\end{equation}
where $|.|^2$ is the usual tensor norm
induced by $g$.

A four-dimensional  KT manifold admits an almost strong KT
structure if and only if it admits a
strong KT. Indeed, (\ref{four1}) for $n=2$ gives
 $\lambda^\Omega=-2\dd\theta \Omega$. In addition, in four dimensions
\begin{equation}\label{four}
T=-*\theta=J\theta\wedge \Omega\ .
\end{equation}
Hence if $\lambda^\Omega=0$, then $\dd\theta=0$ which in turn implies
that $dT=0$.

\section{Curvature Identities}

In this section we shall establish various
 identities for the  curvatures of the Levi-Civita, Chern
and Bismut connections.
In particular we have the following:

\begin{pro}\label{p1}
Let $(M,g,J,\nabla)$  be a KT manifold. The following identities hold
\begin{equation}\label{5}
Ric^g(X,Y) = Ric(X,Y) + \frac{1}{2}\dd T(X,Y) + \frac{1}{4}\sum _{i=1}^{2n}
g\left(T(X,e_i),T(Y,e_i)\right),
\end{equation}
\begin{equation}\label{4}
\rho(X,Y) = Ric (X,JY) +(\nabla_X\theta)JY 
+\frac{1}{4} \lambda^\Omega(X,Y).
\end{equation}
\begin{equation}\label{snov}
b= Scal^{\nabla} - 3\dd\theta -2|\theta|^2 + \frac{1}{3}|T|^2
\end{equation}
\end{pro}

{\it Proof:} Since the torsion is a three-form, we have
\begin{equation}\label{sof}
(\nabla^g_XT)(Y,Z,U) = (\nabla_XT)(Y,Z,U) + \frac{1}{2}
{\sigma \atop XYZ}
\left\{g(T(X,Y),T(Z,U)\right\}.
\end{equation}
Here and henceforth ${\sigma \atop XYZ}$
denote the cyclic sum of $X,Y,Z$.

Using (\ref{1}), we can express  the
 curvature $R^g$ of the Levi-Civita connection in terms of
that of the Bismut connection $R$ as follows:
\begin{eqnarray}\label{15}
R^g(X,Y,Z,U) &=& R(X,Y,Z,U) - \frac{1}{2} (\nabla_XT)(Y,Z,U)
+\frac{1}{2} (\nabla_YT)(X,Z,U)\nonumber \\
             &-&
\frac{1}{2}g(T(X,Y),T(Z,U))
\nonumber\\
&  - & \frac{1}{4}g(T(Y,Z),T(X,U)) -
\frac{1}{4}g(T(Z,X),T(Y,U))\ .
\end{eqnarray}
Taking the trace of  (\ref{15}), using (\ref{sof}) and the fact
that $T$ is a three-form,
we get (\ref{5}).

Further, the exterior derivative $dT$ of $T$ is
given in terms of $\nabla$ by
\begin{eqnarray}\label{13}
dT(X,Y,Z,U)&=&
{\sigma \atop XYZ}\left\{(\nabla_XT)(Y,Z,U)
+2 g(T(X,Y),T(Z,U)\right\}\\
           &-& (\nabla_UT)(X,Y,Z) . \nonumber
\end{eqnarray}
The first Bianchi identity for $\nabla$
together with (\ref{13}) yields
\bea\label{14}
{\sigma \atop XYZ}R(X,Y,Z,U)= dT(X,Y,Z,U)& + & (\nabla_UT)(X,Y,Z)
\nonumber\\
 & - &
{\sigma \atop XYZ}\left\{g(T(X,Y),T(Z,U)\right\}.
\eea
Next we take the trace  of  (\ref{14}) and of  (\ref{13})
with $J$  taking into
account (\ref{3}). Then the
 first equation is  multiplied with two and it is added to the
  second yielding
\bea\label{6}
4\rho(X,Y) & + & 2Ric(Y,JX) -2Ric(X,JY)
\nonumber \\
& = & \lambda^\Omega(X,Y) + 2(\nabla_X\theta)JY -
2(\nabla_Y\theta)JX.
\eea
In addition taking the trace of  (\ref{14}) with $J$ , and using
 (\ref{3})
 and  $R\circ J=J\circ R$, we obtain
\begin{eqnarray}\label{7}
Ric(Y,JX) + Ric(X,JY) &=& \sum_{i=1}^{2n}(R(X,Je_i,e_i,Y)- R(e_i,Y,X,Je_i))\\
&=& -(\nabla_X\theta)JY - (\nabla_Y\theta)JX.\nonumber
\end{eqnarray}
The equation (\ref{4}) in the proposition
 follows from (\ref{6}) and (\ref{7}).

Finally, the equation (\ref{snov}) is
a consequence of (\ref{4}) and the following identity
\begin{equation}\label{eq1}
\sum_{i=1}^{2n}\lambda^\Omega(e_i,Je_i)= 8|\theta|^2 + 8\dd\theta -
\frac{4}{3}|T|^2.
\end{equation}
shown in \cite{AI1}. We remark that the above equation (\ref{eq1})
can be derived by taking the trace of (\ref{13}) twice with $J$.
\hfill {\bf Q.E.D.}

It is straightforward using proposition  \ref{p1} to demonstrate the
following:
\begin{co}
On a KT manifold the Ricci tensor and the Ricci form satisfy the following
relations
\begin{equation}\label{8}
Ric(X,Y)-Ric(Y,X) = - \dd T(X,Y),
\end{equation}
\begin{equation}\label{9}
Ric(JX,JY) - Ric(Y,X) =- (\nabla_{JX}\theta)JY + (\nabla_Y\theta)X,
\end{equation}
\begin{equation}\label{10}
\rho(JX,JY) - \rho(X,Y) = \dd T(JX,Y) - d^{\nabla}\theta(JX,Y),
\end{equation}
where $d^{\nabla}$ is the exterior differential with respect to $\nabla$
given by $ d^{\nabla}\theta(X,Y) = (\nabla_X\theta)Y-(\nabla_Y\theta)X$.
\end{co}
It is worth pointing out that  the Ricci tensor
 of the Bismut connection is not symmetric in general.
Moreover as consequence of (\ref{8}), the Ricci tensor of a linear
connection with totally skew-symmetric
torsion is symmetric
 if and only if the torsion 3-form is co-closed.

The next proposition is our second technical result.
\begin{pro}\label{prs}
On a Hermitian manifold the following identity holds
\begin{equation}\label{c1}
\kappa(JX,Y) = \rho^{1,1}(JX,Y) + <i_XC,i_YC> - \frac{1}{4}
\lambda^\Omega(JX,Y),
\end{equation}
where $\rho^{1,1}$ is the (1,1)-part of
the Bismut-Ricci form, $(i_XC)(Y,Z):=
C(X,Y,Z)$ and $<,>$ denotes the usual
scalar product on tensors induced by $g$.
\end{pro}
{\it Proof:} The above relation can
be most easily established in  local
 holomorphic coordinates $\{z^{\alpha}\}, \alpha =
1,...,n$.
The torsion of the Chern connection is
$$
C_{\alpha\beta\bar{\gamma}}= i d\Omega_{\alpha\beta\bar{\gamma}}=
\partial_{\alpha}g_{\beta\bar{\gamma}} -
\partial_{\beta}g_{\alpha\bar{\gamma}}
$$
and so
$$
\bar{\partial}\partial\Omega = i\left(\partial_{\bar{\delta}}
C_{\gamma\alpha\bar{\beta}} - \partial_{\bar{\beta}}
C_{\gamma\alpha\bar{\delta}}\right)
dz^{\gamma}\wedge dz^{\bar{\delta}}\wedge
dz^{\alpha}\wedge dz^{\bar{\beta}}.
$$
(We have used the Einstein summation conventions.)
The latter equation, using the properties of the
Chern connection,  can be rewritten as
\begin{equation}\label{ch1}
\partial\bar{\partial}\Omega = -i\left(D_{\gamma}
C_{\bar{\delta}\bar{\beta}\alpha} - D_{\alpha}
C_{\bar{\delta}\bar{\beta}\gamma} + C_{\bar{\delta}\bar{\beta}}^{\bar s}
C_{\gamma\alpha\bar{s}}
\right)dz^{\alpha}\wedge dz^{\bar{\beta}}\wedge dz^{\gamma}\wedge
dz^{\bar{\delta}}.
\end{equation}

Taking the trace of  first Bianchi identity of the Chern connection,
$$
K_{\alpha \bar \beta \gamma \bar \lambda} -
K_{\gamma \bar \beta \alpha \bar \lambda} =
- D_{\bar \beta}C_{\alpha \gamma \bar \lambda}
$$
with $g^{\gamma\bar \lambda}$ and after some
computation, one finds   (see e.g. \cite{Bal})
\begin{equation}\label{ch2}
i\rho^D_{\alpha\bar{\beta}} - i\kappa_{\alpha\bar{\beta}} =
-D_{\bar{\beta}}\theta_{\alpha} + D^{\bar s}C_{\bar{s}\bar{\beta}\alpha}.
\end{equation}

Next take the trace of  (\ref{ch1}) and use  (\ref{3}) to get
\begin{equation}\label{sh1}
i(\partial\bar{\partial}\Omega)_{\alpha \bar \beta \gamma \bar \delta}
g^{\gamma \bar \delta}
 = D^{\bar s}
C_{\bar s \bar{\beta}\alpha} + D_{\alpha}\theta_{\bar \beta}
 + C_{\alpha}^{\hspace{1mm}\bar \mu s}
C_{\bar \beta \bar \mu s}.
\end{equation}
Combining (\ref{ch2}) and (\ref{sh1}), we find
\begin{equation}\label{ch3}
i\rho^D_{\alpha\bar{\beta}} - i\kappa_{\alpha\bar{\beta}} =
-D_{\bar{\beta}}\theta_{\alpha} -D_{\alpha}\theta_{\bar{\beta}} -
 C_{\alpha}^{\hspace{1mm}\bar \mu s}
C_{\bar \beta \bar \mu s} +i
(\partial\bar{\partial}\Omega)_{\alpha \bar \beta \gamma \bar \lambda}
g^{\gamma \bar \lambda}.
\end{equation}
Taking the (1,1)-part of (\ref{ricc1}), we derive
\begin{equation}\label{ch4}
i\rho^D_{\alpha \bar \beta} = i\rho_{\alpha \bar \beta}
 -D_{\alpha}\theta_{\bar{\beta}} - D_{\bar{\beta}}\theta_{\alpha},
\end{equation}
since $d(J\theta)_{\alpha \bar \beta}=i(\partial_{\alpha}
\theta_{\bar \beta}+\partial_{\bar \beta}\theta{\alpha}) = i
(D_{\alpha}\theta_{\bar{\beta}} + D_{\bar{\beta}}\theta_{\alpha}).$
We obtain (\ref{c1}) from (\ref{ch4}) and (\ref{ch3}). \hfill {\bf Q.E.D.}

\section{Vanishing Theorems and $SU(n)$ Holonomy}

Suppose that the  (restricted) holonomy of the Bismut connection
$\nabla$ of a KT manifold is contained in  $SU(n)$, ie ${\rm
hol}(\nabla)\subseteq SU(n)$ for short. This holonomy  condition
imposes certain constraints on the geometry,  cohomology groups
and topology of KT manifolds which are similar to those found in
the context of
 compact Calabi-Yau spaces.
In particular, if  ${\rm hol}(\nabla)\subseteq SU(n)$ of a KT manifold $M$,
then the Bismut Ricci form vanishes,
\begin{equation}\label{s1}
\rho = 0.
\end{equation}
If in addition  $M$ is compact the condition (\ref{s1})
implies that the first Chern class
$c_1(M)=0$,  since $\nabla$ gives rise to a
flat unitary connection on the
canonical line bundle $K$.

Recall that for $m>0$ the $m$-th
plurigenenus of a compact complex manifold $(M,J)$
is defined by $p_m(J) = \dim H^0 (M,{\cal O}(K^m))$.

\begin{th}\label{new1}
Let $(M,g,J,\nabla)$ be a  closed 2n-dimensional (non-K\"ahler) KT
manifold with the trace of the Bismut Ricci form to satisfy $$ b>
-|C|^2+{1\over2}h $$ where $2h=\sum_{i=1}^{2n}\lambda^\Omega(Je_i,
e_i)$. Then the plurigenera $p_m(J)=0, \quad m>0$.

\end{th}
{\it Proof:} To prove the statement we shall apply the Gauduchon's plurigenera
theorem \cite{G3,G4}. For this, it is sufficient to show that
\beq\label{pcon}
\int_M\, u_G \,dV_G >0\ ,
\eeq
where   $u_G$
 and $dV_G$ are the trace of the Chern Ricci form with $J$ (\ref{sca}) and the
volume of $M$ with respect the Gauduchon metric $g_G$ of the given
hermitian structure $(g,J)$, respectively.

Let $ g=e^fg_G$ be a conformal transformation
 relating  $g$ to  the Gauduchon metric $g_G$.
For the
functions $u$ and $u_G$ corresponding to  $g$ and $g_G$,
respectively we have (see e.g.\cite{G4})

\begin{equation}\label{tri2}
2 e^f u = 2u_G + n(n-1)<\theta_G,df>_G + n\triangle _G f,
\end{equation}
where $\triangle _Gf$ is the usual Laplacian of $f$ and all terms in the
right hand side are taken with respect to the Gauduchon metric $g_G$.

Taking  the trace in (\ref{c1}), we  get \beq\label{carf}
2u=b+|C|^2-{1\over2}h\ . \eeq Substituting  (\ref{tri2}) into
(\ref{pcon}), integrating over $M$ and
 using $\dd\theta_G=0$ and (\ref{carf}),  we obtain
$$
\int_M u_G \,dV_G ={1\over2} \int_Me^f(b + |C|^2-{1\over2}h)
\,dV_G >0
$$ which is positive under the hypothesis of the
theorem.
 \hfill {\bf Q.E.D.}

An immediate corollary to the above theorem is the following:
\begin{co}\label{costring}

The plurigenera  $p_m(J)=0, m>0$  if any one of the following
conditions holds

\item{(i)} $M$ is an almost strong (non K\"ahler) KT manifold
with ${\rm hol}(\nabla)\subseteq SU(n)$;

\item{(ii)}
 ${\rm hol}(\nabla)\subseteq SU(n)$  and $|C|^2-{1\over2}h > 0$.

\end{co}
{\it Proof:} If the (restricted) holonomy of the Bismut connection
is in $SU(n)$, then $b=0$. In addition under the assumption of
(i), $h=0$ and the corollary follows from the theorem above. Case
(ii) is
 straightforward.
 \hfill {\bf Q.E.D.}

We remark that case (i) of the corollary is applicable to the type
II strings since the relevant manifolds there are strong KT with
(restricted) holonomy of the Bismut connection contained  $SU(n)$.
The case (ii) of the corollary is applicable to heterotic strings.
For heterotic strings the \lq sigma model anomaly cancellation' 
requires that
\beq\label{noncl} 
dH=dT=\mu\big(p_1(TM)-p_1(E)\big)\ , 
\eeq
where $p_1(TM)$ and
$p_1(E)$ are  the Pontriangin classes of the tangent and a vector
bundle $E$ over $M$ in some constant $\mu$
 which depends on the  string tension; the vector bundle $E$
 is associated with the gauge sector of the heterotic string. 
 Both these classes are taken with
respect to Einstein-Hilbert connections; For these connections the
associated curvature is a (1,1) form with respect to the complex
structure $J$ of $M$ and  its trace with $J$ vanishes. In fact
the connection $\tilde \nabla$ on $TM$ is given by setting $T$ to
$-T$ in the definition (\ref{1}) of the Bismut  connection.
Observe that $R(X,Y,Z,W)=\tilde R(Z,W,X,Y)+\frac{1}{2}dT(X,Y,Z,W)$
and note that at the zeroth loop order $dT=0$ . Therefore $$
h={\mu\over4}\left(|\tilde R|^2-|F|^2\right)\ , $$ where $F$ is the Einstein
Hilbert curvature of the vector bundle $E$.

Other applications of the various curvature identities
of the previous section are the following Bochner- Kodaira
vanishing theorem:

\begin{th}\label{new2}
Let $(M,g,J,\nabla)$ be a compact 2n-dimensional KT manifold with
non-negative quadratic form 
$$ <<X,Y>>= \rho^{1,1}(JX,Y) +
<i_XC,i_YC> - \frac{1}{4} \lambda^\Omega(JX,Y)\ . $$
 Then

i) every holomorphic (p,0)-form, $p=1,2,\ldots,n$ is parallel with
respect to the Chern connection;

ii) if moreover $<<,>>$ is positive definite at only one
point, then the Dolbeault cohomology groups $H^0(M,\Lambda^p)=0,
p=1,2,\ldots,n$.

\end{th}
{\it Proof:} To show this, we  apply the vanishing theorem for
holomorphic (p,0)-forms on
 compact Hermitian manifold \cite{KW,Ga}. According to
this general result, it is
sufficient to show that the \lq mean curvature'
 $\kappa$ is non-negative. The non-negativity of $\kappa$ follows from
Proposition~\ref{prs}, formula (\ref{c1}) and the hypothesis
 of the theorem.

\hfill {\bf Q.E.D.}

Two immediate corollaries to the  above theorem are the following:

\begin{co}\label{corn1}
Let $(M,g,J,\nabla)$ be a compact almost strong 2n-dimensional KT manifold
with non-negative (1,1) part $\rho^{1,1}$ of the Bismut Ricci form $\rho$.
Then:

i) every holomorphic (p,0)-form, $p=1,2,\ldots,n$ is parallel
with respect to the Chern connection;

ii) if moreover $\rho^{1,1}$ is positive definite at only one
point, then the Dolbeault cohomology groups $H^0(M,\Lambda^p)=0,
p=1,2,\ldots,n$.

\end{co}

\begin{co}\label{corn2}
Let $(M,g,J,\nabla)$ be a compact 2n-dimensional KT manifold
equipped with a  Bismut connection which  has (restricted)
holonomy contained in  $SU(n)$ and the quadratic form $$ <<X,Y>>=
<i_XC,i_YC> - \frac{1}{4} \lambda^\Omega(JX,Y)\ , $$ is
non-negative. Then every holomorphic (p,0)-form, $p=1,2,\ldots,n$
is parallel with respect to the Chern connection.
\end{co}

The proof of the above two corollaries follow from that of the
main theorem. The first corollary applies to the case of type II
string theory while the second applies to heterotic string theory.
In the former case, if in addition the (restricted) holonomy of
the Bismut connection is contained in $SU(n)$, then existence of
holomorphic (p,0) forms depends on whether at some point in $M$
$<i_XC,i_YC>$ is strictly positive. In the heterotic string case,
it also depends on the positivity properties of $\lambda^\Omega$
which at \lq one loop' is $$
\lambda^\Omega(X,Y)=\mu\big(\sum^{2n}_{i=1}
\big(p_1(TM)-p_1(E)\big)(X,Y,e_i, J e_i)\big)\ . $$

\hfill {\bf Q.E.D..}

Another application of the Theorems ~\ref{new1}
 and Theorem~\ref{new2} above is in context HKT manifolds. For such
spaces the holonomy of the Bismut connection is contained
in $Sp(k)\subset SU(2k)$, $(n=2k)$. In particular
we have the following:
\begin{th}\label{new3}
Let $(M,g,J_a,a=1,2,3, \nabla)$ be a compact almost strong 4n-dimensional
HKT which is not hyperK\"ahler. Then

i) the plurigenera of a complex structure $J_a,a=1,2,3$ is
$p_m(J_a)=0, \quad m>0$;

ii) every holomorphic with respect to a complex structure $J_a,a=1,2,3$
(p,0)-form , $p=1,2,\ldots,2n$ is parallel
with respect to the corresponding Chern connection of $(g,J_a)$
\end{th}

We note that on a compact HKT $p_m(J_a) \in \{0,1\}, a=1,2,3$ as it is
shown in \cite{AI1}.

\subsection{Balanced Hermitian Manifolds}

{\bf Definition} Balanced Hermitian manifolds are Hermitian
manifolds with co-closed K\"ahler form or equivalently with
vanishing Lee form.

Such manifolds have been intensively studied in \cite{Mi,A_B1,A_B2,A_B3}; in
\cite{Ga} they are called semi-K\"ahler manifolds of special type.
This class of manifolds includes the class of K\"ahler manifolds
but also many important classes of non-K\"ahler manifolds, such
as: complex solvmanifolds, twistor spaces of oriented Riemannian
4-manifolds, 1-dimensional families of K\"ahler manifolds (see
\cite{Mi}), some compact Hermitian manifolds with flat Chern
connection (see \cite{Ga}), twistor spaces of quaternionic
K\"ahler manifolds \cite{P,AGI}, manifolds obtained as
modification of compact K\"ahler manifolds \cite{A_B1} and of
compact balanced manifolds \cite{A_B2} (see also \cite{A_B3}).
Some vanishing theorems for balance manifold are given in
\cite{GI2,GI1}.

An application of the previous vanishing theorems is in the
context of balanced hermitian manifolds. For this consider
a hermitian manifold $(M,g,J)$ for which the Lee form
$\theta$ is exact; such manifolds solve the second
Killing spinor equation in (\ref{kse}) and have been investigated
in \cite{SIGP}. If in addition $M$ is compact, then
a direct consequence of the Gauduchon theorem is
that $(M,g_G,J)$ is balanced, where $g_G$ is the Gauduchon metric.
For the existence of almost
strong KT structures on a balanced hermitian manifolds the following
corollary holds:

\begin{co}\label{bala1}
Let $(M,g,J)$ be a 2n-dimensional compact  hermitian manifold for
which the ${\rm hol}(\nabla)\subseteq SU(n)$. If in addition the Lee form
of $(M,g,J)$ is exact, and so $(M,g_G,J)$ is balanced, then the
complex manifold $(M,J)$ does not admit any almost strong KT
structure with ${\rm hol}(\nabla)\subseteq SU(n)$ . In particular if such
a structure do exist, then it is K\"ahler and $(M,J)$ is a
Calabi-Yau space.
\end{co}

{\it Proof:}  This is a direct consequence
of a statement shown in \cite{SIGP}.
 Under the same assumptions as those in  the hypothesis
of the corollary, it has been shown in \cite{strom} (see also \cite{SIGP}),
that $(M,J)$ admits   a globally defined
holomorphic $(n,0)$ form $\tilde \epsilon$ and so $p_1=h^{n,0}\geq 1$.
This form can be expressed in terms of the $\nabla$-parallel
(n,0)-form  $\epsilon$ as
$$
\tilde\epsilon=e^{-f}\epsilon
$$
where $\theta=df$.

Next suppose that there exists a almost strong KT structure on the
complex manifold $(M,J)$ satisfying the conditions of the
corollary. A direct application of corollary~\ref{costring}
reveals that $p_1=0$. Therefore  no such form can exist unless the
torsion of the Chern connection vanishes and the manifold is
Calabi-Yau. \hfill{\bf Q.E.D}

\subsection{Four-dimensional KT Manifolds}

In four dimensions, compact KT manifolds equipped with a Bismut
connection with (restricted) holonomy contained in $SU(2)$ can be
classified. In particular we shall show that such manifolds are
conformal either to a Calabi-Yau space or to a Hopf surface. This
generalizes the result of \cite{rocekb}.

We begin by assuming  that the (restricted) holonomy of the Bismut
connection of a KT manifold $M$ is contained in SU(2)=Sp(1).
Locally this is equivalent to the existence of a HKT structure.
Indeed, there exist (locally) additional two $\nabla $-parallel
almost complex
 structures which together with $J$ satisfy
the relations of imaginary quaternions. The torsion $T$ is
(1,2)+(2,1) form with  respect to each of these new almost complex
structures because the (3,0) and (0,3) parts of a three form on a
four-dimensional manifold vanish identically. This implies that
both new almost complex structures are integrable since their
associated Nijenhuis tensor vanishes. Therefore any KT manifold
with (restricted) holonomy contained in $SU(2)$ admits a (local)
HKT structure.

If the Lee form $d\theta=0$ then the HKT structure is locally
conformally equivalent to a hyper K\"ahler. There are examples of
local HKT structures with non closed Lee form (see e.g.
\cite{Mor,P1}).

On the other hand, the existence of a  HKT structure on a four
manifold is equivalent to the existence of a hypercomplex one
\cite{GT}. This can be seen as follows: Given a  hypercomplex
manifold $(M,J_r)$ equipped with a Riemannian metric $g$, we can
find a  metric $h$ on $M$ which is trihermitian by averaging over
the complex structures, i.e. $$ h(X,Y)=g(X,Y)+\sum_{r=1}^3 g(J_r
X,J_rY)\ . $$ It can be easily seen that $h$ is a  Riemannian
metric on $M$ with the desirable property. Using the integrability
of the complex structures and the fact that (3,0) and (0,3) forms
vanish identically in four dimensions, one can show that
$(M,h,J_r)$ is an HKT manifold.

It is a  well known consequence of the integrability theorem in
\cite{AHS} that the self-dual part of the  Weyl tensor of a
four-dimensional Riemannian manifold $(M,g)$ must vanish if there
exists  a (local) hypercomplex structure on $M$ (see e.g.
\cite{GT} and references there). The converse is not true. For
example, the complex projective space taken with the reverse
orientation is anti-self-dual and does not admit any local
hypercomplex structure \cite{etod}, see also \cite{apgod} and
references there. More precisely, we have

\begin{pro}\label{npa1}
A four-dimensional Hermitian manifold $(M,g,J)$
admits a  anti-self-dual Weyl tensor if and only if
the Bismut Ricci form is anti-self-dual.
\end{pro}
{\it Proof:} A two form on a hermitian surface is anti-self-dual
if it is of type (1,1) and trace-free. Applying (\ref{four}) to
(\ref{10}) we compute $\rho^{(2,0)+(0,2)}=-d\theta_+$, where the
subscript (+) denotes the self-dual part. Hence, the Bismut Ricci
form is of type (1,1) iff the Lee form $\theta$ is anti-self-dual.
It is shown in \cite{AI1} that the trace $b=k$ where $k$ is the
conformal scalar curvature determined by the trace of the
self-dual Weyl tensor $W_+$, ie $k=<3W^+(\Omega),\Omega>$.
Proposition~1 in \cite{Bo2} tells us that $W_+=0$ on a hermitian
surface iff $k=0$  and $d\theta$ is anti-self-dual. \hfill {\bf
Q.E.D.}

In the compact four-dimensional case the situation is different
since there are compact complex surfaces with local hypercomplex
structures which does not admit any global one. This can be seen
in the example of Hopf surfaces explain below.

A Hopf surface $HS$ is by definition a compact complex surface
whose universal covering is ${\bC}^2 -\{0\}$. It was shown by
Kodaira \cite{Kod} that the fundamental group $\pi_1 \tilde{=}
{\bZ}\otimes {\bZ}_n$ and the second Betti number $b_2=0$. A
primari Hopf surface is a Hopf surface with  $\pi_1 \tilde{=}
{\bZ}$. Primary Hopf surfaces are diffeomorphic to $S^1\times
S^3$. Every Hopf surface is finitely covered by a primary one. A
subclass of primary Hopf surface are those  Hopf surfaces $HS_o$
which admit a hermitian metric locally conformally equivalent to a
flat K\"ahler metric. The fundamental group of a primary locally
conformally flat Hopf surface
 is generated by
\cite{God} $\Gamma:(z,w) \rightarrow (az,b^2aw), a,b\in {\cal C},
|b|=1<|a|$. Such a Hopf surface has one hypercomplex structure
exactly when $b=\pm 1$ or when $ab \in \bR$ and it admits exactly
two hyperhermitian structures when both $a$ and $b$ are real
\cite{God}. In view of the classification of compact
hyperhermitian complex surfaces \cite{Bo} these are the only
compact non Calabi-Yau HKT structures.

\begin{th}\label{ps1}

Let $(M^4,g,J)$ be a compact Hermitian surface. The restricted
holonomy of the Bismut connection is contained in SU(2) if and
only if $(M,J)$ is a Calabi-Yau or $(M,g_G,J)$ with the Gauduchon
metric $g_G$ is locally isometric (up to homothety) to $\bR\times
S^{3}$, the Lee form $\theta$ is $\nabla^{g}$-parallel, the
hermitian structure $(g_G,J)$ is locally conformally K\"ahler
flat, and the complex surface $(M,J)$ is a Hopf surface.
 \end{th}
{\it Proof:} If the first Betti number $b_1(M)$ is even then the
surface is of K\"ahler type with zero first Chern class and
therefore it is a Calabi-Yau. Suppose  $b_1(M)$ is odd and
therefore there is no K\"ahler structure on the surface. Using
(\ref{four}), (\ref{four1}) and $\rho=0$ we obtain from (\ref{c1})
that $\kappa = \frac{1}{2}(|\theta|^2+\dd \theta)\Omega$. Hence,
the surface is Einstein-Hermitian and the main result of \cite{GI}
implies the assertion. \hfill {\bf Q.E.D.}

Theorem~\ref{ps1} can be derived also from the properties of Einstein-Weyl 
structures found in \cite{AI2} applying the  connection between
SU(2) holonomy of the Bismut connection with the canonical Einstein-Weyl 
structure on a hermitian surface discovered in \cite{AI1}.

The description of the Hopf surfaces and Theorem~\ref{ps1} show
that there exist compact hermitian surfaces with restricted 
 holonomy of the
Bismut connection contained in SU(2) which does not admit a HKT
structure.

We finish this section with the following non-vanishing result which 
gives obstructions to the existence of non Calabi-Yau spaces with 
SU(n) holonomy.

\begin{th}\label{thnov}
Let $(M,J)$ be a 2n-dimensional compact complex manifold with vanishing 
first Chern class and does not admit any K\"ahler metric. If there exist 
a hermitian structure $(g,J)$ such that the restricted holonomy group of the
Bismut connection is contained in SU(n) then either

i) $h^{0,n}=h^{n,0}=1$ if the Lee form is exact

or

ii) $h^{0,1}\ge 1$ if the Lee form is not exact.
\end{th}
{\it Proof.} If the Lee form is an exact 1-form then the Gauduchon metric is 
balanced and $h^{n,0}=1$ by the result in \cite{strom} (see also \cite{SIGP}).

Suppose that the Lee form $\theta$ is not exact ie the Gauduchon metric 
$g_G=e^fg$ is not balanced and the corresponding Lee form $\theta_G \not=0$.
The equation (\ref{ricc1}) shows that the 2-form $dJ\theta$ is an (1,1)-form
since $\rho=0$ and $\rho^D$ is an (1,1)-form. This yields $d\theta$ is an 
(1,1)-form ie $\bar{\partial}\theta^{0,1}= \bar{\partial}(J\theta)^{0,1}=0$,
where $^{0,1}$ means the (0,1)-part. Then, we get 
$\bar{\partial}\theta_G^{0,1}=0$ due to the equality 
$\theta_G^{0,1}=\theta^{0,1}+(n-1)\bar{\partial}f$. Combining the latter with 
$\dd_G \theta_G=0$ we obtain that $\theta_G^{0,1}$ is a $\bar{\partial}$-
harmonic with respect to the Gauduchon metric non zero form. Hence, 
$h^{0,1}\ge 1$. \hfill {\bf Q.E.D.}

\section{String equations and Hermitian Manifolds}

We shall investigate the cohomological properties
of solutions of the string equations
 (\ref{st}) which are
complex  manifolds $(M,g,J)$. For this the torsion $T$
of the Bismut connection $\nabla$ is related
to the three-form field strength $H$  of string theory as
\begin{equation}\label{s0}
H= T\ .
\end{equation}
As it was already mentioned, for applications in string theory the
restricted holonomy of the Bismut connection is contained in
$SU(n)$. In this case the relation between string equations and
complex manifolds can be made manifest.

We shall seek to describe vanishing theorems for solutions of the
string equations which are KT manifolds for which the Bismut
connection has (restricted) holonomy contained in $SU(n)$. Most of
the statements below apply to almost strong KT manifolds.
Therefore most of our results concern type II strings for which
$H$ is closed. Our results apply also to the heterotic string
but the contribution from the gauge sector has been neglected.

The presence of a strong KT structure on
$(M,g,J)$ that is associated
with type II strings allows the use
of theorem shown in \cite{AI1}
which we shall state here without proof as follows:

\begin{th}\label{ASI}
Let $(M,g,J)$ be a compact 2n-dimensional $(n>1)$ strong KT
manifold with K\"ahler form $\Omega$. Suppose that the Lee form
$\theta$ is co-closed, $\dd\theta=0$ and that the (1,1)-part of
the Ricci form of the Bismut connection is non-negative everywhere
on $M$.

a) Then every $\bar\partial$-harmonic $(0,p)$-form, $p=1,\dots,n$,
 is parallel with respect to the Bismut connection.

b) If moreover the (1,1)-part of the Ricci form of the Bismut connection
is strictly positive at some point, then the cohomology groups
$H^p(M,{\cal O})$ vanish for $p=1,\dots,n$.

\end{th}

Note that the above theorem has been  stated in a different but
equivalent way in \cite{AI1}.
The solutions of the string equations
can be separated into two classes depending on whether
the dilaton is a constant or not.

\subsection{Constant Dilaton}

Let the dilation $\phi$ be constant.
In such case the string equations (\ref{st}), after using the
above identification (\ref{s0})
of $H$ and $T$, become
\bea
Ric^g_{ij}-{1\over4} T_{imn} T_j{}^{mn}&=&0
\nonumber\\
\nabla^g_i T^{imn}&=&0\ .
\eea
or equivalently
\beq\label{st1}
Ric=0\ .
\eeq
In particular the above equation implies that $\dd T=0$ (see (\ref{8})).
At present, we take $dT\not=0$. The restrictions imposed on $T$
necessary for  the various theorems will be stated explicitly.

Assuming that the (restricted) holonomy of the Bismut connection
is in $SU(n)$, and so $\rho=0$, we find using (\ref{4}) that the
string equations (\ref{st1}) can be rewritten as
\begin{equation}\label{st1'}
(\nabla_X\theta )Y = \frac{1}{4}\lambda^\Omega(X,JY).
\end{equation}
For manifolds with  an almost strong KT structure
$\lambda^\Omega=0$ and so the above equation implies  that the Lee
form $\theta$ is parallel with respect to the Bismut connection,
i.e.
$$ \nabla \theta=0\ .
$$
In particular,
$\theta^{\#}$ is a Killing vector field. This can be easily seen
using the definition of the Bismut connection and  (\ref{1a}).

{\bf Remark 2.} For almost strong KT manifolds,
$\lambda^\Omega=0$,  the equation (\ref{eq1}) shows that
 the Lee form is identically zero iff
$(M,g,J)$ is a K\"ahler manifold.
Hence, on any non-K\"ahler almost strong KT
 manifold which is a solution of
(\ref{st1}),
there is a globally defined non-zero Killing
vector field $\theta^{\#}$.

There is an alternative way to characterize the string equations
the following:

\begin{th}\label{th1}
Let $(M,g,J,\nabla)$ be a 2n-dimensional compact strong KT manifold and
${\rm hol}(\nabla)\subseteq  SU(n)$.
Then $(M,g,J,\nabla)$ is a solution of the string equations
(\ref{st1}) if and only if  the scalar curvature
of the Bismut connection vanishes, $Scal^{\nabla}=0$.
\end{th}
{\it Proof:}{}From the assumptions of the theorem we have
$$
\rho=dT=0\ .
$$
 Since $dT=0$, $\lambda^\Omega=0$. Substituting
this into equation (\ref{eq1}), we find
\beq\label{eqq1}
2|\dd\theta|^2+2|\theta|^2 -{1\over3} |T|^2=0\ .
\eeq

To show the theorem in one direction, we
insert $\rho=Ric=\lambda^\Omega=0$   into (\ref{4}) and  find that
$$
\nabla \theta =0
$$
which in turn  implies $\dd\theta=0$; therefore $g$
is the Gauduchon metric.
Substituting this into (\ref{eqq1}), we get
\beq\label{qqq1}
2|\theta|^2 -{1\over3} |T|^2=0\ .
\eeq
Moreover since  $\rho=0$, then  $b=0$. Substituting this, (\ref{qqq1})
and $\dd\theta=0$   into  (\ref{snov}),
we find
$Scal^\nabla=0$. Which proves the theorem in one direction.

For the converse, suppose that $Scal^{\nabla}=0$. Then, substituting
$b=0$ and (\ref{eqq1}) into
 (\ref{snov}) yields
\beq\label{coc} \dd\theta=0\ . \eeq Next taking the trace of
(\ref{4}) using $\rho=Scal^{\nabla}=\lambda^\Omega=0$, we find
that $\sum_{i=1}^{2n} (\nabla_{e_i} \theta) Je_i=0$. In turn using
(\ref{3}) this implies that \beq\label{coca} \sum (\nabla^g_{e_i}
\theta) Je_i=0\ . \eeq From (\ref{coc}) and (\ref{coca}), we
conclude that $\theta^{0,1}$ is $\bar\partial$ co-closed, ie $$
\bar\partial^\dagger \theta^{0,1}=0\ , $$ where the superscript
${}^{0,1}$ means the (0,1)-part of the form. Note that the Lee
form $\theta$ is not identically zero because of Remark 2. Since
$\rho=0$, equation (\ref{ricc1}) implies that $\rho^D=d(J\theta)$
and therefore $d(J\theta)$ is an (1,1)-form. This in turn implies
that $d\theta$ is a (1,1) from, ie $$
\bar{\partial}\theta^{0,1}=0\ . $$ Thus $\theta^{0,1}$ is both
$\bar\partial$ closed and co-closed and therefore $\bar\partial$
harmonic. Applying  the vanishing theorem of \cite{AI1} stated in
\ref{ASI}, we conclude that $$ \nabla \theta=0\ . $$ Substituting
this together with $\rho=\lambda^\Omega=0$ into (\ref{4}), we get
$Ric=0$. Hence, we recover  the string equations for constant
dilation $\phi$.

\hfill {\bf Q.E.D.}

We remark that  the solutions of the string equations that we
investigate are {\it not} required to satisfy the second Killing
spinor equation in (\ref{kse}). If they did, then the second in
(\ref{kse}) implies that $\theta=2d\phi$ (see \cite{strom, SIGP}).
Since the dilaton $\phi$ is constant, this in turn implies that
$\theta=0$. If $(M,g,J)$ satisfies the assumptions of the theorem
above, substituting $\theta=0$ in  (\ref{eqq1}), we find that the
torsion $T$ vanishes and the manifold $M$ is Calabi-Yau.

There are some restrictions on the  cohomology of
manifolds that are solutions of
the string equations (\ref{st1}).
More precisely the following result holds:
\begin{th}\label{thh1}
Let $(M,g,J,\nabla)$ be a 2n-dimensional compact strong KT manifold and
be a solution of the string equations (\ref{st1}). Moreover let
${\rm hol}(\nabla)\subseteq  SU(n)$. Then  $(M,g,J,\nabla)$
has the following properties:

i) Every $\bar{\partial}$-harmonic (0,p)-form,
$p=1,2,..., n$ is parallel with
respect to the Bismut connection $\nabla$. Therefore
$h^{0,p}=dimH^p(M,{\cal O}) \le {n \choose p}$, $h^{0,1}\ge 1$ and
the dimension of the
space of Killing vector fields is at least 2$h^{0,1}\ge 2$;

ii) every holomorphic (p,0)-form , $p=1,2,\ldots,2n$ is parallel
with respect to the Chern connection and therefore
$h^{p,0}= \le {n \choose p}$;

iii) the plurigenera $p_m(J)=0, \quad m>0$ provided $(g,J)$ is not
K\"ahler.
\end{th}

{\it Proof:} To show (i), we apply the vanishing
theorem of \cite{AI1} stated in \ref{ASI} using
the assumptions  that $\rho=dT=0$.

Moreover using the assumptions $\rho=\lambda^\Omega=b=0$, (ii)
and (iii) are direct consequence of Theorems~\ref{new2} and
~\ref{new1} , respectively.
\hfill {\bf Q.E.D.}

\begin{co}
Let $(M,J,g)$ be as in  Theorem~\ref{thh1}.  If in addition
the Lee form  $\theta$ is closed,  then
$\theta$ is $\nabla^g$-parallel
and $(M,g)$ is locally isometric to
$N^{2n-1}\times \bR$, where
$N^{2n-1}$ is a (2n-1)-dimensional Riemannian
manifold with non-negative Riemannian Ricci curvature.
\end{co}
{\it Proof:}Using  the assumptions of the theorem, we have shown that
$\theta$ is parallel with respect to the Bismut connection.
This and the additional assumption that $\theta$ is closed imply
that
$$
i_{\theta^{\#}}T=0\ .
$$
Substituting this back into $\nabla \theta=0$, we find
that $\theta$ is parallel with respect to the Levi-Civita
connection. The non-negativity of the Riemannian Ricci
curvature follows from (\ref{5}) and the assumptions
of the corollary.
\hfill {\bf Q.E.D.}

There are various applications to the above theorems. One such
application is in the context of HKT manifolds. Since the
 holonomy of the Bismut connection for HKT manifolds
is contained in $Sp(k)$, $n=2k$, it is also contained in $SU(n)$.
In particular as a direct application of Theorem~\ref{th1}, we
have the following:

\begin{th}\label{th2}
A compact strong HKT manifold is a solution of the string
equations with constant dilaton (\ref{st1}),
iff $Scal^{\nabla}=0$.
In addition the statements (i), (ii) and (iii) of
Theorem~\ref{thh1} taken
with respect to any of the complex
structure $J_a, a=1,2,3$ also hold.
\end{th}

As another application, we can relate KT manifolds
and HKT manifolds as follows:

\begin{co}\label{co1}
Let $(M,g,J,\nabla)$ be a 4k-dimensional compact strong KT
manifold, the (restricted)  holonomy of $\nabla$ be in SU(2k) with
$Scal^{\nabla}=0$. If there exists a non-degenerate
$\bar{\partial}$-harmonic (0,2)-form $\phi$, then $(M,g,J,\nabla)$
is a strong HKT manifold.
\end{co}

{\it Proof:} If there is such a (0,2)-form $\phi$, then
 applying  Theorem~\ref{thh1} we conclude that it is parallel
with respect to the Bismut connection. Hence, the (restricted)
holonomy of $\nabla$ is contained in Sp(k). \hfill {\bf Q.E.D.}

We remark that the existence of $\bar{\partial}$-harmonic (0,2)-form
is not necessary for the existence of a $\nabla$-parallel one.
For example,  $SU(3)$ with any left invariant complex
 structure has $h^{0,2}=0$
but there exist a left invariant HKT structure on this group. The
same example can be used to demonstrate that $h^{4,0}=p_1(J)= {\rm
dim} H^0\big(M, {\cal O}(K)\big)=0$ despite the fact that there is
a  $\nabla$-parallel $(4,0)$-form because the 
holonomy of $\nabla$ is trivial in accordance with
Theorem~\ref{thh1}.

{\bf Examples.} Compact strong solutions are even dimensional
Lie groups endowed with an left invariant
complex structure compatible with a bi-invariant
Riemannian metric \cite{SSTV,OP1}
The induced Bismut connection is flat and the torsion is parallel
and so both  closed and
co-closed. The Lee form is not always closed in these examples.
Some of these group manifold examples admit in fact
a  HKT structure \cite{SSTV,OP1,GrP}.

In four dimensions,  (\ref{four}) implies
$\nabla^g\theta =\nabla \theta$.
The conditions $dT=0, \dd T=0$ are
equivalent to $\dd\theta =0, d\theta=0$,
respectively. The standard hermitian structure
on the Hopf surfaces $HS_o$ has
$\nabla^g$-parallel Lee form and $\rho =0$.
Then $Ric=0$ by (\ref{4}) and therefore it is a
solution of the string equations (\ref{st1}). In fact,
the Hopf surfaces  $HS_o$ with the standard hermitian structure
 are the unique compact non
Calabi-Yau solution of the string equation with constant dilation
and (restricted) holonomy of the Bismut connection contained in
SU(2) by Theorem~\ref{ps1}. In particular, every compact strong
HKT surface solves the string equation (\ref{st1}) with constant
dilation.

\subsection{Non-constant dilaton}

Let $(M,g,J)$ be a KT manifold equipped with a Bismut connection
with (restricted) holonomy contained in $SU(n)$. In addition, we
take  the dilation $\phi$ not to be a constant, $\phi\not={\rm
const}$. The string equations (\ref{st}) can be written using
(\ref{s0}) and (\ref{5}) as
\begin{eqnarray}\label{st2}
& &Ric + \frac{1}{2}\dd T + 2\nabla^g d\phi=0,\\
& &\dd T + 2i_{d\phi^\#}T=0,\nonumber
\end{eqnarray}
where  $i_X$ denotes
the interior multiplication by a vector $X$.

Using (\ref{4}),  (\ref{10}) and the assumption that the
(restricted) holonomy of the Bismut connection is contained in
$SU(n)$, so we have $\rho=0$,  we find  that the string equations
(\ref{st2}) can be expressed in terms of the Lee form $\theta$  as
\begin{eqnarray}\label{st2'}
-(\nabla_X\theta)Y - (\nabla_Y\theta)X +\frac{1}{2}\lambda^\Omega(X,JY)
 + 4(\nabla^g_Xd\phi)Y& = &0,\\
(\nabla_X\theta)Y - (\nabla_Y\theta)X + 2T(d\phi^\#,X,Y)& = &0.\nonumber
\end{eqnarray}

The first string equation in (\ref{st2'}) can be also written in the 
following equivalent way 
\begin{equation}\label{cnew}
(\nabla_X\eta)Y + (\nabla_Y\eta)X =
\frac{1}{2}\lambda^\Omega(X,JY),
\end{equation}
where the 1-form $\eta$ is given by
\begin{equation}\label{stef1}
\eta = \theta - 2d\phi.
\end{equation}

{\bf Remark 3.} For  compact non K\"ahler solutions of
(\ref{cnew}) the Lee form
$\theta $ cannot be identically zero. Indeed, combining (\ref{cnew})
 together with
(\ref{eq1}), one finds
$$ 6\dd\theta - \frac{2}{3}|T|^2+4|\theta|^2 -
 4\dd d\phi=0\ .
$$
Integrating this equality over a compact manifold without boundary,
we obtain $T=0$ if and only if  $\theta=0$.

\vskip 0.2truecm

We have
\begin{th}\label{noco}
A 2n-dimensional KT manifold $(M,g,J)$ is a (local) solution to
the first string equation of (\ref{st}) and to the Killing spinor
equations (\ref{kse}) with non-constant dilation if and only if
$(M,g,J)$ is an almost strong KT manifold with closed Lee form and
${\rm hol}(\nabla)\subseteq SU(n)$.  If in addition $M$ is an eight-dimensional
 compact space, then it admits a   strong KT structure.
\end{th}
{\it Proof:} If the Killing spinor equations (\ref{kse}) are
satisfied, then the second in (\ref{kse}) implies that
$\theta=2d\phi$ (see \cite{strom, SIGP}). {}From (\ref{stef1}), we find
that $\eta=0$ and the conditions of the theorem imply (\ref{cnew}).

The last statement in the theorem follows because in eight dimensions
the condition
$\lambda^{\Omega}=0$ and the fact that $dT$ is a (2,2) form imply that
 $dT$ is
self-dual 4-form with respect to the Hodge * operator. Now if $M$ is
compact,  we get $\int_M|dT|^2 \,dV=-\int_M<*d*dT,T>
\,dV=-\int_M<*d^2T,T> \,dV=0$. So $M$ is a strong KT manifold.
 \hfill {\bf Q.E.D.}

We remark that any manifold that 
satisfies the assumptions of the above theorem
and  equation (\ref{cnew}) admits an almost
strong KT structure. This appears to be in conflict with
the heterotic five-brane solution of \cite{callan}. However
this is not the case because the contribution of
the gauge sector in the string equations (\ref{st}) have been neglected;
although the contribution of the gauge sector due to the
sigma model anomaly  in the non-closure 
of $T$ (\ref{noncl}) has been
taken into account. In fact it has been shown in \cite{howepap}
that for the heterotic five-brane to solve the 
string equations, the two loop contribution 
should be taken into account. Such corrections to the
string equations
involve quadratic terms in the curvatures and they have not
be taken into account in this paper.

Next using (\ref{10}) together with $\rho=0$,  we find that
\begin{equation}\label{ns1}
\dd T=d^{\nabla}\theta = d\theta - i_{\theta^{\#}}T.
\end{equation}
Combining (\ref{ns1}) with the second string equation in
(\ref{st2}), we get $$ d\theta = i_{(\theta^{\#} - 2d\phi^{\#})}T.
$$ Finally, we observe that  if $(M,g,J)$ is a KT manifold and the
(restricted) holonomy of $\nabla$  is contained in SU(n), then the
second string equation in (\ref{st}) is equivalent to the equation
\begin{equation}\label{ster}
(\nabla_X\eta)Y -
(\nabla_Y\eta)X  = 0.
\end{equation}
Combining (\ref{cnew}) with (\ref{ster}) we get that the string 
equations (\ref{st}) are equivalent to the following one
\begin{equation}\label{stef}
\nabla_X \eta(Y) =\frac{1}{4}\lambda^\Omega(X,JY),
\end{equation}

Thus we have shown the following theorem:

\begin{th}\label{noc}
Let $(M,g,J)$ be a 2n-dimensional KT manifold and ${\rm
hol}(\nabla)\subseteq SU(n)$. $(M,g,J)$ is a solution to the
string equations (\ref{st}) if and only if the equation
(\ref{stef}) holds.
\end{th}

If in addition $(M,g,J)$ is an almost strong KT
manifold again with $\rho=0$, as
it will be in the case of type II strings,
then the string equations
are equivalent to
\beq\label{stst1}
\nabla\eta=0
\eeq
and so $\eta$ is $\nabla$-parallel.
 There are two special cases that one can consider
the following:

{\bf Case 1.} Let us suppose that the Lee form is closed i.e. $$
d\theta =0\ . $$ Then locally $\theta =df$ for a smooth function
$f$. Taking $$ \phi = \frac{1}{2}f $$ we obtain a local solution
of the string equations (\ref{stst1}) for which $\eta=0$. Thus,
every almost strong KT manifold with (restricted) holonomy in
SU(n) and closed Lee form gives a local solution of the string
equations (\ref{st}) and (\ref{kse}) with non-constant dilation
$\phi$.

Now suppose that  in addition the class $[\theta]$ of the Lee form
in $H^1(M)$ is trivial, ie $\theta$ is exact. As it has already
been mentioned this is precisely the case that the second Killing
spinor equation in (\ref{kse}) admits a solution. Therefore only
backgrounds with $\theta$ an exact form are supersymmetric.
 Moreover let $(M,g,J)$ be a 2n-dimensional compact
almost strong KT manifold with ${\rm hol}(\nabla)\subseteq SU(n)$
Then one can apply the corollary ~\ref{bala1} shown in the
context of balanced hermitian manifolds to show that $(M,J)$ is
Calabi-Yau.

\vskip 0.1cm

{\bf Case 2.} Alternatively, the Lee form $\theta$ may not be
closed,  $d\theta \not=0$. If $(M,g,J)$ is an almost strong KT
manifold for which the associated Bismut connection has
(restricted) holonomy contained  in $SU(n)$,
 then as we have already mentioned in (\ref{stst1}) the string equations
are equivalent to the condition that the form $\eta$ is
$\nabla$-parallel . In particular $\eta^{\#}$ is a non-zero
Killing vector field and $$ d\eta=d\theta\ . $$ Conversely, if an
almost strong KT manifold $(M,g,J)$ for which the associated
Bismut connection has (restricted) holonomy contained  in $SU(n)$
admits a $\nabla$-parallel one-form $\eta$ such that
$d\eta=d\theta$, then $(M,g,J)$ is a local solution of the string
equations with dilaton  $\phi =-\frac{1}{2}f$, where $f$ is
determined by $df=\eta-\theta$. Thus we have shown the following
corollary:

\begin{co}
Let $(M,g,J)$ be an almost strong 2n-dimensional KT manifold for
which the associated Bismut connection has (restricted) holonomy
contained in $SU(n)$ and the Lee form is not closed. Then
$(M,g,J)$ is a solution of the string equations with non-constant
dilaton $\phi=(1/2) f$ if and only if it admits a
$\nabla$-parallel one-form $\eta$ such that $d\eta=d\theta$ and
$\eta-\theta=df$ is an exact one-form.
\end{co}

Next suppose that $(M, g,J)$ is a four-dimensional
hermitian manifold.  Using (\ref{four1}), we get
that string equations (\ref{stef}) take the form
\begin{equation}\label{four2}
\nabla \eta =\frac{1}{2}\dd\theta\otimes g,
\end{equation}
 $(M,g,J)$
where $\eta$ is given by (\ref{stef1}). In particular,
$\eta^{\#}$ is a non-zero
conformal Killing vector field provided $\dd\theta \not=0$.

\vskip 0.4truecm

\section{Concluding Remarks}

We have found conditions for the 
existence of  Bismut connections
on hermitian  manifolds (KT) for which their (restricted) holonomy
contained in $SU(n)$. These 
conditions can be expressed in terms
of the vanishing of certain 
cohomology groups. For example under
certain additional assumptions the plurigenera vanish.
 We also consider various applications
of our results  in the context of string theory
and in the context of balanced hermitian manifolds.

Despite the various developments the last few years in the
context of KT and HKT manifolds, 
the existence of an  HKT structure on a
manifold has not been expressed in terms of conditions on its 
cohomology groups in
parallel with similar  developments in the case of
hyper-K\"ahler manifolds. For this  
an analogue of the Calabi-Yau
theorem is needed for the case of
 KT manifolds which will involve
hermitian  manifolds with topologically 
trivial canonical bundle. It is not
clear for example under which 
conditions a KT manifold with trivial canonical
bundle admits a KT structure 
associated with Bismut connection for
which its (restricted) holonomy is contained in $SU(n)$. 
We have found in Theorem~\ref{thnov} 
obstructions to the existence of KT 
structure with holonomy in SU(n) on a 
compact complex manifold with vanishing 
first Chern class of non-K\"ahler 
type.  There seem
though to be counter examples for the existence of strong KT
structures of this type.

An alternative but essentially equivalent way
to state the Calabi-Yau type of conjecture mentioned
above which seems to be of importance
for the further development of the 
theory is the following: Given a 2n-dimensional 
$(n>2)$ compact complex manifold with 
zero first Chern class and $h^{n,0}=1$ 
or $h^{0,1}\ge 1$, does there exist a 
hermitian metric with vanishing Ricci 
form of the Bismut connection? 
The dimension $(n>2)$ in the above 
question is essential since the Inoe 
surface has vanishing first Chern class 
and $h^{0,1}=1$ but it does not 
admit Hermitian structure with SU(2)
 holonomy of the Bismut connection 
by Theorem~\ref{ps1}.

\vskip 0.4truecm

{\bf Acknowledgments:} A part of this research was done during the
visit of (S.I.) at the Abdus Salam International Centre for
Theoretical Physics, Trieste Italy, Spring 2000 and 
Humboldt University, Berlin, Fall 2000.
S.I. thanks the
Abdus Salam ICTP and HU Berlin for support and the excellent environment.
 S.I is supported by the contracts MM
809/1998 of the Ministry of Science and Education of Bulgaria,
238/1998 of the University of Sofia ``St. Kl. Ohridski'' and EDGE,
Contract HPRN-CT-2000-00101.  G.P. is supported by a University
Research Fellowship from the Royal Society. This work is partially
supported by SPG grant PPA/G/S/1998/00613.

\end{document}